\documentclass[11pt,reqno]{amsart}
\usepackage{amsmath,amsfonts,amssymb,amscd,amsthm,amsbsy}
\newtheorem{thm}{Theorem}
\newtheorem{lem}[thm]{Lemma}

\def\DeltaNTA{\text{\rm $\Delta$-N.T.A.}}
\def\div{\operatorname{div}}
\def\Im{\operatorname{Im}}
\def\Lip{\operatorname{Lip}}
\def\tw{\widetilde w}
\def\ep{\varepsilon}
\def\L{{\mathcal L}}
\def\KK{{\mathbb K}}
\def\nat{{\mathbb N}}
\def\real{{\mathbb R}}
\def\negint{\mathop{\int\mkern-19mu {\diagup}}\nolimits}
\def\smallnegint{\mathop{\int\mkern-13mu
\raise.5ex\hbox{${\scriptscriptstyle\diagup}$}}\nolimits}
\begin{document} 
\title[Structure of the free boundary for lower dimensional obstacle problems]
{The structure of the free boundary for lower dimensional obstacle problems}
\author[I. ATHANASOPOULOS, L. A. CAFFARELLI, S. SALSA]
{I. Athanasopoulos,$^1$ L. A. Caffarelli,$^2$ S. Salsa$^3$} 
\address{I. Athanasopoulos\\
University of Crete\\
Department of Applied Mathematics\\
Knossos Avenue\\
71 409 Herakleion Crete, Greece {\rm and}
Institute of Applied and Computational Mathematics\\
Forth, P.O. Box 1385\\
71110 Herakleion Crete, Greece}
\email{athan@tem.uoc.gr}
\address{L.A. Caffarelli\\
Department of Mathematics\\
University of Texas at Austin\\
Austin, TX~78712} 
\email{caffarel@math.utexas.edu}
\address{S. Salsa\\
Dipartimento di Matematica\\
Politecnico di Milano\\
Piazza Leonardo da Vinci, 32\\
20133 Milano, Italia}
\email{sansal@mate.polimi.it}
\thanks{${}^1$ Supported by RTN Project HPRN-CT-2002-00274.}
\thanks{${}^2$ Supported by National Science Foundation Grant DMS-0140338.}
\thanks{${}^3$ Supported by Murst. Polimi, 2004.}
\begin{abstract}
We study the regularity of the ``free surface'' in boundary obstacle 
problems. We show that near a non-degenerate point the free boundary 
is a $C^{1,\alpha}$ $(n-2)$-dimensional surface in $\real^{n-1}$.
\end{abstract}
\maketitle
\section{Introduction}

The purpose of this paper is to study the structure and regularity of the 
free boundary in ``boundary obstacle'' problems.
Boundary obstacle problems concern the following question.
\begin{quote}
We are given a smooth $\Omega$ in $\real^n$, $n\ge 3$, and seek a 
function $u$ that:  
\begin{itemize}
\item[a)] In the interior of $\Omega$, $u$ satisfies a nice, elliptic 
equation, say $\Delta u =f$.
\item[b)] Along the boundary of $\Omega$, instead of giving Dirichlet or 
Neumann conditions we prescribe ``complementary conditions'' of the 
following type. 
As long as $u$ is bigger than some prescribed function $\varphi$, 
there is no flux across $\partial \Omega$: $u_\nu =0$.
But as soon as $u$ becomes equal to $\varphi$, boundary flux, $u_\nu$,  
is turned on $(u_\nu >0)$ to keep $u$ above $\varphi$.
\end{itemize}
\end{quote}

This type of problem arises in elasticity (the Signorini problem) when 
an elastic body is at rest, partially laying on a surface, 
in optimal control of temperature across a surface 
(see \cite{F}, \cite{A}), in the modelling 
of semipermeable membranes  where some saline  concentration can flow 
through the membrane only in one direction
(see Duvaut-Lions \cite{DL}) 
and in financial math when the random variation of underlying asset 
changes in a discontinuous fashion 
(a Levi process) (see \cite{S} and references there).

There is considerable literature on the regularity properties of the 
solution  (see \cite{F}, \cite{C}, \cite{R}, \cite{U}). 
In particular, two of the authors proved recently (see \cite{AC2}) 
the optimal regularity of solutions to such a problem.

This opens the way to study the properties of the interface
by using geometric P.D.E.\ techniques.
This is precisely what we develop in this paper. 
We show that there is one basic global non-degenerate profile after 
blow up, and that in a neighborhood of a point that has this profile 
the free boundary is a $C^{1,\alpha}$ ``curve'' on the boundary 
(i.e., an $n-2$ dimensional graph on the $n-1$  dimensional boundary).

Simple examples show that singular free boundary points and degenerate
profiles are unavoidable. 
For simplicity, in this paper we only treat the case in which $\partial\Omega$
is locally a hyperplane and $f,\varphi\equiv 0$.

\section{Description of the problem and known results}

In this section we explain exactly which kind of problem we shall deal 
with and we recall some known results.

Let $B_1 = B_1(0)$ the unit ball in $\real^n$, $n\ge 2$; we write 
points $x\in \real^n$ as $x= (x',x_n) \in \real^{n-1}\times \real$ 
and denote by $\Pi$ the hyperplane $\{(x',x_n) :x_n=0\}$. 

Given a smooth function $\varphi$ on $\partial B_1$ we look at the unique 
minimizer $u$ of the Dirichlet integral 
$$J(v) = \int_{B_1} |\nabla v|^2$$
over the closed convex set 
$$\KK = \{v\in H^1 (B_1),\ v= \varphi\text{ on } \partial B_1,\ 
v(x',0) \ge 0\}$$
The minimizer  $u$ can be constructed also as the least 
superharmonic function in $\KK$.
To have a nontrivial coincidence set $\Lambda(u) = \{(x',0) :u(x',0)=0\}$
and a nontrivial free boundary $F(u)$, the boundary of 
the set $(u\ge \varphi)$ on 
$\Pi\cap B_1$, we assume that $\varphi$ changes sign and that 
$\varphi (\theta',0) >0$, $\theta' \in\partial B'_1 = \Pi \cap \partial B_1$.
Without losing generality we can choose $\varphi$ symmetric with respect 
to the hyperplane $\Pi$ so that, $u$ also is  symmetric with respect to $\Pi$
(otherwise we can symmetrize without changing the coincidence set). 

The solution $u$ is harmonic in $B_1\setminus \Lambda(u)$, it is globally 
Lipschitz continuous, and (see \cite{C}) 
\begin{equation}\label{eq2.1}
\|u\|_{\Lip (B_{1/2})} \le C\|u\|_{L^2(B_1)}
\end{equation}
Moreover,
\begin{equation}\label{eq2.2} 
\inf_{B_{1/2}}u_{\tau\tau} \ge -C\|u\|_{L^2(B_1)}
\end{equation}
for every direction $\tau$ on $\Pi$. 
We will call {\em tangential\/} such directions.
Inequality \eqref{eq2.2} expresses semiconvexity of $u$ along 
tangential directions. 

The optimal regularity of $u$, proven in \cite{AC2}, is $C^{1,1/2}$ on 
either side of $\Pi$, and 
\begin{equation}\label{eq2.3} 
\|u\|_{C^{1,1/2} (B_{1/2}^\pm)} \le C\|u\|_{L^2(B_1)}\ .
\end{equation}
Furthermore, $u_{x_n} =0$ on $\{(x',0): u(x',0)>0\}$ and 
$u_{x_n} (x',0+) \le 0$ on $\Lambda (u)$.

In this paper we want to examine the structure of the free boundary 
$F(u)$ (clearly in dimension $n\ge3$) through the analysis of 
asymptotic profiles around one of its points, that we assume to be
the origin.

It turns out that only in correspondence to a specific asymptotic profile 
(that we call nondegenerate) it is possible to achieve smoothness of $F(u)$.
To get a clue of what happens let us start with an observation of 
Hans Lewy in dimension~2. 

The complex function $w = u_x -u_y$ $(x_n=y)$ is analytic outside 
$\Lambda (u)$, thus 
$$u_x u_y = \Im (w^2)$$
is harmonic  and vanishes on $y=0$.
Thus $u_x u_y$ has a harmonic odd extension across $y=0$ and $w^2$ 
has an analytic extension. 
Then $w$ is $C^{1/2}$ and $u\in C^{1,1/2}$, which is indeed the optimal 
regularity.
Accordingly, the first admissible nontrivial global solution is 
$u_0(x) = \rho^{3/2} \cos \frac32  \theta$ and this is the typical 
nondegenerate asymptotic profile. 
On the other hand there are solutions like 
$\rho^{k +1/2}\cos ((k+1/2)\theta)$, 
$k \in\nat$, $k >1$, or $\rho^{2k} \cos 2k\theta$,
$k \ge1$, with higher  order asymptotic behavior. 

In correspondence to points with these asymptotic profiles the free 
boundary could be very narrow or a singular point. 
Notice that these 2-dimensional solutions can be considered as 
$n$-dimensional solutions, constant with respect to the other 
$n-2$ variables, so that analogous considerations can be made 
in any dimension.

\section{Monotonicity Formulas}

In this section we prove some monotonicity formulas that play a crucial 
role in the identification of limiting blow-up profiles.

\begin{lem}[Almgreen's frequency formula]\label{lem:frequency}
Let $u$ be a  continuous function on $\bar B_r$, harmonic in $B_r^+$, 
$u(0)=0$, $u(x',0) \cdot u_{x_n} (x',0)=0$.
Define, for $0<r<1$,
$$D_r (u) = r \frac{\int_{B_r} |\nabla u|^2}{\int_{\partial B_r} u^2 d\sigma}
 \equiv r \frac{V_r}{S_r}\ .$$
Then, for $0<r\le \frac12$, $D'_r(u) \ge 0$ $({}' = \frac{d}{dr})$.
Moreover, let 
$$\mu = \lim_{r\to 0^+} D_r (u)\ .$$
then $D'_r(u) \equiv 0$ in $(0,\frac12)$  if and only if 
$$u(x) = |x|^\mu g(\theta)\qquad \theta \in \partial B_1$$
and $\mu \ge \frac32$.
\end{lem}

\begin{proof}
We have 
$$\log D_r = \log r + \log V_r - \log S_r$$
and 
$$\frac{d}{dr} \log D_r = \frac1r + \frac{V'_r}{V_r} - \frac{S'_r}{S_r}\ .$$
By rescaling, it is enough to show that 
\begin{equation}\label{eq1}
1 + \frac{V'_1}{V_1} - \frac{S'_1}{S_1} \ge 0
\end{equation}
that is ($\nu$ exterior normal)
\begin{equation}\label{eq2}
2- n + \frac{\int_{\partial B_1} |\nabla u|^2\,d\sigma}
{\int_{B_1} |\nabla u|^2} - 2
\frac{\int_{\partial B_1} u\, u_\nu\, d\sigma}
{\int_{\partial B-1} u^2\, d\sigma}\ . 
\end{equation}

Since $u(x',0) = u_{x_n}(x',0) =0$ we get, after an integration by parts, 
$$\int_{\partial B_1} |\nabla u|^2 = \frac12 
\int_{\partial B_1} \Delta (u^2) = \int_{\partial B_1} u\, u_\nu\,d\sigma\ .$$
To control $\int_{\partial B_1} |\nabla u|^2\,d\sigma$ we use the 
divergence theorem in $B_1\setminus \Lambda(u)$. 
Let 
$$h(x) = \div [ x|\nabla u|^2 - 2(x\cdot \nabla u) \nabla u]\ .$$
Notice that, in our case 
$$h(x) = (n-2) |\nabla u|^2\ .$$
{From} Gauss formula, we have (using that on $\Lambda$ $u_\tau$ vanishes 
continuously)
\begin{equation}\label{eq3}
(n-2) \int_{B_1} |\nabla u|^2 = \int_{B_1} h 
= \int_{\partial B_1} |\nabla u|^2\, d\sigma 
- 2 \int_{\partial B_1} u_\nu^2 \,d\sigma\ .
\end{equation}
By inserting \eqref{eq3} into \eqref{eq2} we obtain 
$$1 + \frac{V'_1}{V_1} - \frac{S'_1}{S_1}
= 2 \frac{\int_{\partial B_1} u_\nu^2 \,d\sigma}
{\int_{\partial B_1} u\, u_\nu \, d\sigma} 
- 2 \frac{\int_{\partial B_1} u\, u_\nu\,d\sigma}
{\int_{\partial B_1} u^2\, d\sigma} 
\ge 0$$
by Schwarz inequality. 
The equality sign in $D'_r (u) = 0$ holds for $0<r\le \frac12$
if and only if $u$ is proportional to $u_\nu$ on 
$\partial B_r$ for every $r$, which implies $u$ is of the form 
$$u(x) = h(|x|) g(\theta) \qquad \theta \in \partial B_1\ .$$
{From} the radial formula of the Laplace operator, in a 
neighborhood of any point where $u\ne 0$, it must be 
$$h(|x|) = |x|^\mu\ .$$
In fact, by unique continuation, $\mu$ must be the same 
for all components of $B_1\setminus \Lambda(u)$
where $g$ has constant sign.
Thus, each connected component of the region where 
$u$ is harmonic is a cone, generated by the support of $g$.
Finally, from optimal regularity, it must be $\mu \ge \frac32$.
\end{proof}

An important consequence is the following result.

\begin{lem}\label{lem2}
Let $u$ and $\mu$ as in Lemma~\ref{lem:frequency}.
Let 
$$\varphi (r) = \negint_{\partial B_r} u^2\,d\sigma\qquad 
0<r\le 1\ .$$

{\rm (a)} $r^{-2\mu} \varphi (r)$ is increasing and 
$[r^{-2\mu} \varphi (r)]' = 0$ in $(0,1)$
if and only if 
$$u(x) = |x|^\mu g(\theta)\qquad\theta \in \partial B_1$$
with $\mu \ge \frac32$.
\bigskip

{\rm (b)} Let $0<r<R\le 1$; given $\ep >0$, for $r\le r_0(\ep)$
\begin{equation}\label{eq4}
\varphi (R) \le \left( \frac{R}{r}\right)^{2(\mu +\ep)} \varphi (r)\  .
\end{equation}
\end{lem}

\begin{proof}
(a) We have
\begin{equation}\label{eq5} 
\varphi'(r) = \frac{d}{dr} \negint_{\partial B_r} u^2 
= 2\negint_{\partial B_r} u\, u_\nu\,d\sigma 
= 2r \negint_{B_r} |\nabla u|^2
\end{equation}
so that 
$$\frac{d}{dr} [r^{-2\mu} \varphi (r)] 
= 2r^{-2\mu -n} \bigg\{ r\int_{B_r} |\nabla u|^2 
- \mu \int_{\partial B_r} u^2\,d\sigma\bigg\}$$
and (a) follows from the frequency formula.

(b) Let $r_0 = r_0(\ep)$ such that $D_r(u) \le \mu +\ep$.
{From} \eqref{eq5}
$$D_r (u) = \frac{r}2\, \frac{d}{dr} \log \varphi (r) \le \mu +\ep$$
and \eqref{eq4} follows by integrating over $(r,R)$.
\end{proof}

\section{Limiting Profiles}

Given a solution $u$ of our thin obstacle problem, 
we consider the blow-up family
$$v_r (x) = \frac{u(rx)}{(\smallnegint_{\partial B_r} u^2)^{1/2}}\ .$$
If $\mu = \lim_{r\to 0^+} D_r(u)$, our purpose is to identify the 
limit of $v_r$ as $r\to 0$ when $\frac32\le \mu \le 2$.

Observe that 
\begin{equation}\label{eq6}
\|v_r\|_{L^2 (\partial B_1)} = 1
\end{equation}
and, from Lemma~\ref{lem2}
$$\|v_r\|_{L^2(B_R)} \le R^{(\mu +\ep)}$$
for every $R>1$ and every small $r$ such that $rR \le r_0(\ep)$.
Thus, a sequence $v_j = v_{r_j}$ converges in $L^2$ and 
uniformly on every compact set in $\real^n$ to a nontrivial 
(because of \eqref{eq6}) global solution $v_0$.

Since
$$D_{r_j} (u) = D_1 (v_{r_j})\to D_1 (v_0) = \mu$$
as $r_j \to 0$, from Lemma~\ref{lem:frequency} we deduce that 
$$v_0 (x) = |x|^\mu g(\theta) \qquad \theta \in \partial B_1\ .$$
We now distinguish two cases.

\subsection*{The case $\frac32 \le \mu <2$}
{From} the tangential quasi-convexity property of $u$, we have, 
for every tangential direction $\tau$:
\begin{equation}\label{eq7}
D_{\tau\tau} v_{r_j} \ge 
- c\, \frac{r_j^2}{(\smallnegint_{\partial B_{r_j}} u^2)^{1/2}}\ .
\end{equation}
In Lemma~\ref{lem2}(b), choose $\ep$ such that $\mu +\ep <2$.
Then, letting $r_j\to 0$ in \eqref{eq7} we obtain from \eqref{eq4} 
$$D_{\tau\tau} v_0 \ge 0$$
so that $v_0$ is tangentially convex and $\Lambda(v_0)$ is a 
convex cone. 
We first observe that on $\Lambda_0$, $v_0\equiv 0$, and $D_nv_0\le0$ 
and for $x_n \ge 0$, $D_{nn}\le 0$. 
This implies that $v_0 (x',x_n)\le0$ if $(x',0)$ belongs to $\Lambda_0$. 
Assume now that the vector $-e_{n-1}$ belongs to $(\Lambda_0)^0$. 
For any point $x$, consider the line $L_x = \{x+te_{n-1}\}$. 
For $t$ negative enough the function $v_0 (x+te_{n-1})$ becomes negative 
from the remark above. 

Since $v_0$ is convex along $L_x$, it follows that $w= D_{e_{n-1}}v$ cannot be 
negative anywhere on $L_x$. 
In particular, since $x$ is arbitrary $w\ge0$ in $R^n$. 
On the other hand, $w=0$ on $\Lambda (v_0)$ and $w_{x_n} =0$ on 
${\{x_n = 0\}\setminus \Lambda (v_0)}$ (by symmetry). 
Thus, the restriction of $w$ to the unit sphere 
must be the first eigenfunction of the Dirichlet problem for 
the spherical Laplacian, with zero data on $\partial B_1 \cap \Lambda (v_0)$. 
Now, if $\Lambda (v_0)$ is not a half-plane, and thus 
from convexity, is strictly contained in half a plane. 
Then the homogeneity degree of $w$ should be less than $1/2$ 
(see \cite{AC2}), since homogeneity $1/2$ corresponds to the case in which
half a plane is removed, contradicting $\mu \ge \frac32$. 
Therefore $\Lambda (v_0)$ is a half-plane, 
$w(x) = \rho^{1/2}\sin \frac{\psi}2$ where ${\rho^2 = x_{n-1}^2 + x_n^2}$, 
and $\tan \psi = x_n/x_{n-1}$. 
This  implies 
\begin{equation}\label{eq3.8} 
v_0 (x) = \rho^{3/2} \cos \frac32 \psi\ .
\end{equation}
Observe that if $\tau = \alpha e_{n-1} + \beta e$, where $e$ is tangential, 
$e\perp e_{n-1}$ and $\alpha^2 + \beta^2 =1$, $\alpha >0$, then outside 
a $\eta$-strip, $|x_n| <\eta$, we have 
\begin{equation}\label{eq3.9} 
D_\tau v_0 (x) \ge C(\alpha) \eta^{1/2}\ .
\end{equation}

\subsection*{The case $\mu =2$}
The limiting profile is of the form $v_0 (x) = |x|^2 g(\theta)$, 
$\theta \in\partial B_1$ and $\Lambda (v_0)$ is a cone. 
Consider $w= D_{x_n} v_0$; $w$ is linearly homogeneous and $w=0$ on 
$\{x_n =0\}\setminus \Lambda (v_0)$. 
We reflect evenly with respect to the hyperplane $x_n =0$, defining 
\begin{equation*}
\widetilde w (x) = \begin{cases}
w(x',x_n)&x_n >0\\
\noalign{\vskip6pt}
w(x',-x_n)&x_n <0 
\end{cases}\ .
\end{equation*}
Suppose $\tw$ changes sign. 
Then, since $\tw$ is harmonic on its  support and $w(0)=0$, we can apply the  
monotonicity formula in \cite[Theorem 12.3]{CS}, to $w^+$ and $w^-$. 
According to this formula, the linear behavior of $\tw$ forces $\tw$ 
to be a two plane solution with respect to a direction transversal to 
the plane $x_n=0$, say, $\tw (x) = \alpha x_{n-1}^+ - \beta x_{n-1}^-$, 
due to the even symmetry of $\tw$. 
This is a contradiction  since along $x_n=0$, $w$ is negative on $\Lambda_0$ 
and zero otherwise
and therefore $\tw$ cannot change sign. 
Suppose now that $\Lambda (v_0)$ has non-empty interior. 
Then $\tw$ is the first eigenfunction for the spherical Laplacian, with zero 
boundary data on $(\{x_n =0\}\setminus \Lambda (v_0)) \cap \partial B_1$. 
This forces a superlinear behavior of $\tw$ at the origin 
since linear behavior corresponds to a half sphere  and we reach 
again a contradiction. 
Thus, $\Lambda (v_0)$ has empty interior, $v_0$ is harmonic across 
$\Lambda (v_0)$ and therefore $v_0$ must coincide with a quadratic 
polynomial ($v_0 (x) = \sum_{i<n} a_i x_i^2 - C  x_n^2$, $a_i \ge 0$).

\section{Lipschitz continuity of the Free Boundary ($\mu <2$)}

Through the identification of the limiting profile in section~4, we can 
prove that, when $\frac32 \le \mu <2$, the free boundary $F(u)$ is locally 
a Lipschitz graph.
Precisely:

\begin{lem}\label{lem3}
Let $u$ be a solution of the thin obstacle problem in $B_1$. 
Assume that  $\frac32 \le \mu <2$. 
Then, there exists  a neighborhood of the origin $B_\rho$ and 
a cone of tangential directions $\Gamma'$ $(e_{n-1},\theta)$, 
with axis $e_{n-1}$  and opening $\theta \ge \frac{\pi}3$ (say), such that, 
for every $\tau\in\Gamma'$ $(e_{n-1},\theta)$, we have 
$$D_\tau u\ge 0\ .$$
In particular, in that neighborhood, $F(u)$ is the graph of a Lipschitz 
function $x_{n-1} = f(x_{n-1},\ldots,x_{n-2})$.
\end{lem}

\begin{proof} 
As in Lemma~\ref{lem:frequency}, let 
$$v_{r_j}(x) = \frac{ u(r_jx)}{\smallnegint_{\partial B_{r_j}}(u^2)^{1/2}}\ .$$
We know from section 3 that $v_{r_j} (x) \to v_0(x)$ with $v_0$ given by 
\eqref{eq3.8}, uniformly on compact sets. 
Fix $\alpha \ge \frac12$ (say) and let $\tau = \alpha e_{n-1} +\beta e$ 
be a tangential direction $(\alpha^2 +\beta^2 =1)$. 
For $\sigma >0$, small, and $r_j\le r_0(\sigma)$, we deduce from \eqref{eq3.9}
that $D_\tau v_{r_j}$ enjoys the following properties in $B_{5/6}$.
\begin{itemize}
\item[(i)] $D_\tau v_{r_j} \ge 0$ outside the strip $|x_n| <\sigma$; 
\item[(ii)] $D_\tau v_{r_j} \ge c_0 >0$ for $|x_n| \ge \frac12$;
\item[(iii)] $D_\tau v_{r_j} \ge -c\sigma^{1/2}$ in the strip $|x| <\sigma$
\end{itemize}
(from optimal regularity).

Then, we conclude the proof by applying to $h= D_\tau v_{r_j}$ the following 
approximation Lemma. 
\end{proof}

\begin{lem}\label{lem:approx}
Let $u$ be a solution of the thin obstacle problem in $B_1$. 
Suppose $h$ is a continuous function with the following properties:
\begin{itemize}
\item[(i)]  $\Delta h \le 0$ in $B_1\setminus \Lambda (u)$;
\item[(ii)] $h\ge 0$ for $|x_n| \ge \sigma$, $h=0$ on $\Lambda (u)$, 
with $\sigma >0$, small; 
\item[(iii)] $h\ge c_0 >0$ for $|x_n| \ge \frac1{8(n-1)}$;
\item[(iv)] $h > -\omega (\sigma)$, where $\omega$ is the modulus of 
continuity of $h$, for $|x_n| <\sigma$. 
\end{itemize} 
There exists $\sigma_0 = \sigma_0$ $(n,c_0,\omega)$ such that, if 
$\sigma \le \sigma_0$ then $h\ge0$ in $B_{1/2}$.
\end{lem}

\begin{proof} 
Suppose $z = (z',z_n) \in B_{1/2}$ and $h(z) <0$. 
Let 
$$Q = \Big\{ (x',x_n) : |x'-z'| \le \tfrac13\ ,\ |x_n| < 
\frac1{4(n-1)} \Big\}$$
and 
$$P(x',x_n) = |x'-z'|^2 - (n-1) x_n^2\ .$$
Define
$$v(x)= h(x) + \delta P(x)$$
where $\delta >0$ is to be chosen later. 
We have 
\begin{itemize}
\item[(a)] $v(z) = h(z) - \delta (n-1) z_n^2 <0$
\item[(b)] $\Delta v\le 0$\quad outside $\Lambda (u)$
\item[(c)] $v\ge 0$ on $\Lambda (u)$, since $h\ge 0$, 
$P\ge 0$ there.
\end{itemize}
Thus, $v$ must have a negative minimum on $\partial Q$.

On $\partial Q \cap \{ |x_n| > 1/8(n-1)\}$, 
$$v\ge c_0 - \frac{\delta}{16 (n-1)} \ge 0$$
if $\delta \le 16 (n-1) c_0$.

On $|x'-z'| =1/3$,  $\sigma \le |x_n| \le 1/8(n-1)$, we have 
$h\ge 0$ so that 
$$v\ge \delta \left[ \frac19 - \frac1{64(n-1)} \right] \ge 0\ .$$
Finally, on $|z'-z'| = 1/3$, $|x_n| <\sigma$, we have 
$$v\ge c\omega (\sigma) + \delta \left[ \frac19 - (n-1) \sigma^2
\right] \ge 0$$
if $\sigma$ is small enough.

Hence, $v\ge 0$ on $\partial Q$ and we have reached a 
contradiction. 
Therefore $h\ge 0$ in $B_{1/2}$.
\end{proof}

\section{Boundary Harnack Principles and the $C^{1,\alpha}$ 
Regularity of the Free Boundary $(\mu <2)$}

We are now in position to show that the free boundary is locally 
a $C^{1,\alpha}$ graph, if $\mu <2$.
Precisely, our main result is the following.

\begin{thm}\label{thm:main}
Let $u$ be a solution of the thin obstacle problem in $B_1$. 
If ${\mu = \lim_{r\to 0+} D_r (u) <2}$ 
then the free boundary $F(u)$ is given in a neighborhood of 
the origin by the graph of a $C^{1,\alpha}$ function 
$x_{n-1} = f(x_1,\ldots, x_{n-2})$.
\end{thm}

One way to prove the theorem is to use the results in \cite{AC1}.
Through a bilipschitz transformation, a neighborhood of the origin
in $B_1\setminus \Lambda (u)$ is mapped onto the upper half ball,
say, $B^+ = \{ |z| <1$, $z_{n-1} >0\}$, and the Laplace operator is 
transformed into a uniformly elliptic divergence form operator. 
Each tangential derivative $D_\tau u$, with $\tau$ belonging to 
the cone $\Gamma' (e_{n-1},\theta)$ of monotone directions, 
is mapped onto a positive solution of $\L v=0$ in $B^+$, vanishing 
on $\{z_{n-1} =0\}$.
An application of Corollary~1 in \cite{AC1} concludes the proof.

On the other hand, there is a more direct proof based on the 
following result, that could be of interest in itself.

Let $D$ be a subdomain of $B_1$ and let $\Omega = \partial D
\cap B_1$. 
We denote by $d_g(x,y)$ the geodesic distance in $D$ of the 
points $x,y$.
We will assume that the following properties hold:
\begin{itemize}
\item[(1)] For every $x,y\in D$, $d_g(x,y)$ is finite.
\item[(2)] {\bf Non tangential ball condition}. 
Let $Q\in \Omega$. 
There exist positive numbers $r_0 = r_0 (D,Q)$ and $\eta =\eta (D)$
such that, for every $r\le r_0$ there is a point $A_r (Q) \in B_r(Q)$ 
such that 
$$B_{\eta r} (A_r (Q)) \subset B_r (Q) \cap D\ .$$
\item[(3)] {\bf Harnack chain condition}. 
There exists a constant $M = M(D)$ such that, for all $x,y\in D$,
$\ep >0$ and $k\in \nat$ satisfying 
$$d(x,\Omega) >\ep\ ,\quad d(y,\Omega) > \ep\ ,\quad 
d_g  (x,y) < 2^k \ep\ ,$$
there is a sequence of $Mk$ balls $B_{r_1},\ldots,B_{r_{Mk}} \subset D$
with 
$$x\in B_{r_1}, \quad y\in B_{r_{Mk}}, \quad B_{r_j} \cap B_{r_{j+1}}
\ne \emptyset \qquad (j=1,\ldots, Mk-1)$$
and 
$$\frac12 r_j < d (B_{r_j} , \Omega) < 4 r_j\qquad 
(j=1,\ldots,Mk)\ .$$
\item[(4)] {\bf Uniform capacity condition}. 
Let $Q\in\Omega$. 
There exist positive numbers $r_0 = r_0 (D,Q)$ and $\gamma = \gamma (D)$
such that 
$$\text{cap}_\Delta \left[ (B_{2^{-k}r} (Q) \setminus 
B_{2^{-k-1}r} (Q) ) \cap \Omega\right] \ge \gamma r^{n-2}$$
where $\text{cap}_\Delta (k)$ is the capacity of $k$ in $B_1$, 
with respect to the Laplace operator.
\end{itemize}
\medskip

\noindent
Conditions (2) and (3) appear in the notion of more tangentially 
accessible domain (see \cite{JK}). 
Condition~(4) replaces the exterior tangential ball property in that 
definition.
Since condition~(4) is related to the Laplace operator we call a 
domain $D$ with properties (1)--(4) a $\DeltaNTA$\ domain.
A simple example of $\Delta$-N.T.A.\ domain is an 
$(n-1)$-dimensional smooth manifold with Lipschitz boundary.

Let now $\L$ be a uniformly elliptic operator with ellipticity constant 
$\lambda$ and bounded measurable coefficients. 
Recall that $\text{cap}_\L (K) \sim \text{cap}_\Delta (K)$ with 
constant depending only on $\lambda$ and $n$, so that the 
notion of $\DeltaNTA$\ domains is actually related to an 
entire class of operators. 
The following result holds.

\begin{thm}\label{thm:bdryharnack}
{\rm (Boundary Harnack Principles).} 
Let $D\subset B_1$ be a $\DeltaNTA$\ domain. 
Supposes $v$ and $w$ are positive functions in $D$, 
continuously vanishing on $\Omega$, satisfying $\L v= \L w=0$ 
in $D$. 
Assume $x_0 \in D\cap B_{2/3}$, $d(x_0,\Omega) = d_0 >0$ and 
$v(x_0) = w(x_0)=1$. 
Then:
\begin{itemize}
\item[(a)] For every $Q\in\Omega \cap B_{1/2}$
$$\sup_{D\cap B_{1/3}(Q)} v \le C(u,d,d_0,D)\qquad 
\text{\em (Carleson estimate)}$$
and 
$$\sup_{D\cap B_{1/3}(Q)} \frac{v}{w} \le C(u,\lambda,d_0,D)\ .$$
\item[(b)] $\frac{v}w$ is H\"older continuous in $B_{1/2} \cap D$
up to $\Omega$.
\end{itemize}
\end{thm}

\begin{proof}
the proof follows by now standard lines (see for instance 
\cite[section 11.2]{CS} and \cite{JK}).
We sketch the main steps emphasizing the main differences.
\medskip

\noindent (a)
Fix $Q\in\Omega\cap B_{1/2}$ and let 
$$v(y_0) = N= \sup_{B_{1/3}Q\cap D} v\ .$$
The interior ball condition and the Harnack chain condition plus the interior 
Harnack inequality imply that if $N$ is large, $d(y_0,\Omega) \equiv 
|y_0 - Q_0| \le N^{-\ep}$ where $\ep = \ep (n,\lambda,d_0,D)>0$. 
Let $r_0  = d(y_0,\Omega)$. 

The uniform capacity condition implies that 
$$\sup_{B_{2r_0}(Q_0)} v \equiv v (y_1) \ge CN$$
where $C = C(n,\lambda,D) >1$. 
Iterating the process, one constructs a sequence of points 
$y_k$, satisfying 
\begin{itemize}
\item[i)] $v(y_k) \ge C^k N$
\item[ii)] $d(y_k,\Omega) \le (C^k N)^{-\ep}$
\item[iii)] $|y_k - y_{k+1}| \le 4 (C^k N)^{-\ep}$.
\end{itemize}
If $N$ is large enough, we can make 
$$\sum |y_k - y_{k+1}| \le \frac1{16}$$ 
and we get a contradiction. 
This proves (a).
To prove (b),  let $P\in B_{1/3}(Q)\cap \Omega$ and 
$R_0 = d(x_0,P) + \frac{d_0}2$. 
Notice that $\frac32 d_0 \le R_0 \le 1$. 
Define 
$$\psi_{R_0} (P) = B_{R_0} (P) \cap D$$ 
and 
$$\Sigma_0 = \partial B_{R_0} (P) \cap B_{d_0} (x_0)\ .$$ 
Observe that $\Sigma_0 \subset D$. 
We first control the Green's function $G(x,x_0)$ for $\L$ 
in $\psi_{R_0} (P)$ from above by the $\L$-harmonic 
measure $\omega_\L^x (\Sigma_0)$, in $\psi_{R_0}(P)
\setminus B_{d_0/3} (x_0)$. 
This follows from the maximum principle.
In fact, on $\partial B_{d_0/3} (x_0)$ we have, from H\"older continuity, 
$$\omega_\L^x (\Sigma_0) \ge c>0$$
and 
$$G(x,x_0)\le c\,d_0^{2-n}\ .$$
On the other hand, on $\partial \psi_{R_0} (P)$, we have 
$G(x,x_0)=0$ and $\omega_\L^x (\Sigma_0) \ge 0$. 
Therefore, outside $B_{d_0/3} (x_0)$ we get 
\begin{equation}\label{star}
G(x,x_0) \le c\, d^{2-n} \omega_\L^x (\Sigma_0)\ .
\end{equation}
Let now $\Sigma_1 = \partial \psi_{R_0} (P) \setminus \Omega$
and let $\varphi$ be a $C^\infty$ cut-off function such that 
$\varphi \equiv 0$ in $B_{R_0/4} (P)$, $\varphi \equiv 1$ outside 
$B_{R_0/2}(P)$ and $0\le \varphi \le 1$ in $B_{R_0/2} (P)
\setminus B_{R_0/4} (P) \equiv C_{R_0} (P)$. 

We have 
$$\omega_\L^x (\Sigma_1) \le \int_{\partial \psi_{R_0}(P)} 
\varphi \,d\omega_\L^x\ .$$
Fix $x\in B_{R_0/8} (P) \cap D$. 
Then 
$$0 = \varphi (x) = \int_{\partial \psi_{R_0}(P)} \varphi \, 
d\omega_\L^x 
- \int_{C_{R_0}(P)\cap D} a_{ij} (x) D_{y_i} G(x,y) 
D_{y_j} \varphi (y)\,dy\ .$$
Therefore, from Caccoppoli estimate and Carleson estimate, 
we have, in $B_{R_0/8} (P)\cap D$, 
\begin{equation}\label{starstar}
\begin{split}
\omega_\L^* (\Sigma_1) 
& =  c(u,\lambda, d_0) \bigg( \int_{C_{R_0}(P)\cap D} 
|\nabla_y G (x,y) |^2\,dy\bigg)^{1/2}\\
\noalign{\vskip6pt}
&\le c(u,\lambda,d_0) \bigg( \int_{C_{R_0}(P)\cap D} 
G^2 (x,y)\, dy\bigg)^{1/2}\\
\noalign{\vskip6pt}
& \le c(u,\lambda,d_0) G (x,x_0)\ .
\end{split}
\end{equation}
{From} \eqref{star} and \eqref{starstar} we obtain the following 
doubling condition for the $\L$-harmonic measure: 
$$\omega_\L^x (\Sigma_1) \le c(u,\lambda,d_0) \omega_\L^x 
(\Sigma_0)$$
for every $x\in B_{R_0/8} (P) \cap D$.

The rest of the proof of (a) and the proof of (b) follow now, for
instance, as in \cite[section 11.2]{CS}.
\end{proof}

\begin{proof}[Proof of Theorem~\ref{thm:main}]
We apply Theorem~\ref{thm:bdryharnack} with $\Omega = \Lambda (u)$, 
$D= B_1\setminus \Lambda (u)$ and $v= D_\tau u$, $w= D_{e_{n-1}}u$
where $\tau \in \Gamma' (e_{n-1},\theta)$. 
We obtain, in particular, that on $\{x_n=0\}\setminus \Lambda (u)$, 
the quotient $D_\tau u/D_{e_{n-1}} u$ is H\"older continuous up 
to $F(u)$ in a neighborhood of the origin. 
This implies that the level sets in $\real^{n-1}$ of $u$ are $C^{1,\alpha}$ 
surfaces and, in particular, 
the $C^{1,\alpha}$ regularity of $F(u)$ in $B_{1/2}$.
\end{proof}


\begin{thebibliography}{MMM}

\bibitem[A]{A}
I. Athanasopoulos, 
{\em Regularity of the Solution of an Evolution Problem with Inequalities 
on the Boundary}, 
Comm. P.D.E. {\bf 7} (1982), 1453--1465.

\bibitem[AC1]{AC1}
I. Athanasopoulos, L. A. Caffarelli, 
{\em A Theorem of Real Analysis and its Applications to Free Boundary Problems},
Comm. in P.A.M.  {\bf 38}  (1985), 499--502.

\bibitem[AC2]{AC2}
\bysame,
{\em Optimal Regularity of Lower Dimensional Obstacle Problems},



\bibitem[C]{C}
L. A. Caffarelli, 
{\em Further regularity for the Signorini problem}, 
Comm. P.D.E. {\bf4}(9) (1979), 1067--1075.


\bibitem[CS]{CS} 
L. A. Caffarelli, S. Salsa, 
{\em A Geometric Approach to Free Boundary Problems}, 
American Mathematical Society, Providence, G.S.M., Vol.68.

\bibitem[DL]{DL}
G. Duvaut, J.L. Lions, 
{\em Les inequations en mechanique et en physique}, 
Paris Dunod 1972.

\bibitem[F]{F} 
J. Frehse, 
{\em On Signorini's problem and variational problems with thin obstacles},
Ann. Scuola Norm. Sup. Pisa {\bf4} (1977),  343--362.

\bibitem[JK]{JK}
D. Jerison, C. Kenig, 
{\em Boundary behavior of harmonic functions in non-tangentially 
accessible domains},
Adv. in Math. {\bf 46}, No.1 (1982), 80--147.

\bibitem[R]{R} 
D. Richardson, Thesis, University of British Columbia, 1978.

\bibitem[S]{S} 
L. Silvestre, Thesis, University of Texas at Austin, 2005.

\bibitem[U]{U}
N.N. Uraltseva, 
{\em On the regularity of solutions of variational inequalities} 
(Russian), Usp. Mat. Nauk {\bf42} (1987), 151--174.


\end{thebibliography}
\end{document}